\documentclass{amsart}
\usepackage{graphicx}
\newtheorem{thm}{Theorem}[section]

\newtheorem{lem}[thm]{Lemma}

\theoremstyle{definition}
\newtheorem{defn}[thm]{Definition}
\theoremstyle{remark}

\numberwithin{equation}{section}

\begin{document}

\title{A Geometric Characterization of Rational Groups}%
\author{Cecil Andrew Ellard}%
\address{Bloomington, Indiana}%
\email{cellard@ivytech.edu}%

\thanks{}%

\begin{abstract}
We give a geometric characterization of finite rational groups. In particular, we prove that a finite group is rational if and only if there exists a finite geometry $\Gamma$ of type $I$ and action of $G$ on $\Gamma$ as a group of automorphisms such that if $g$ and $h$ are elements of $G$ fixing the same number of flags of type $J$ for all subsets $J$ of $I$, then $g$ and $h$ are conjugate in $G$. 

\end{abstract}
\maketitle

\section{Rational Groups}

A finite group $G$ is said to be a $\textit{rational group}$ if every every ordinary (i.e., complex) irreducible character $\chi$ of G is a  rational-valued. Since ordinary character values of a finite group are algebraic integers, and since any rational algebraic integer is an integer, it follows that the character table of a finite rational group consists entirely of integers. Examples of rational groups include the Coxeter groups of type $A_{n},$ (which are finite symmetric groups), $B_{n}, D_{n}, E_{6}, E_{7}, E_{8}, F_{4},$ and $G_{2}$. Since each of these Coxeter groups acts in a natural way on the coset geometry of its maximal parabolic subgroups, it seems natural to seek geometric equivalents to rationality for finite groups.\\

Let $I$ be a nonempty set. An $\textit{incidence geometry}$ (or just a $\textit{geometry}$) of type $I$ is a triple $\Gamma = (X, *, t)$ where $X$ is a non-empty set (the \textit{objects} of the geometry), * is a symmetric and reflexive relation on $X$ (the \textit{incidence} relation) and $t:X \rightarrow I$ is a function (the \textit{type} function) such that if  $t(x)=t(x^{\prime})$ and $x*x^{\prime}$ then $x=x^{\prime}$.  The geometry is said to be finite if $X$ is finite. A $\textit{flag}$ of the geometry is a set $F$ of mutually incident elements of $X$, and the $\textit{type}$ of a flag $F$ is the set of types of the elements of $F$. An automorphism of $\Gamma$ is a bijection $\alpha : X \rightarrow X$ such that $\alpha$ and $\alpha^{-1}$ preserve incidence and types. The group of all automorphisms of $\Gamma$ will be denoted by $Aut(\Gamma)$, and we will say that $G$ \textit{acts on} $\Gamma$ if there is a homomorphism $\theta : G \rightarrow Aut(\Gamma)$. When $G$ acts on $\Gamma$, for $g \in G$ and $J \subseteq I$, we will define $Fix_{J}(g)$ to be the set of flags of type $J$ fixed by $g$. The purpose of this note is to prove the following:

\begin{thm} Let $G$ be a finite group. Then the following are equivalent:\newline
\newline
(1) $G$ is a rational group.\newline
\newline
(2) There exists a finite geometry $\Gamma$ of type $I$ and action of $G$ on $\Gamma$ as a group \newline 
\indent $\hskip 3pt$  of automorphisms such that if $g$ and $h$ are elements of $G$ fixing the same \newline 
\indent $\hskip 3pt$  number of flags of type $J$ for all subsets $J$ of $I$, then $g$ and $h$ are conjugate in $G$. \\

\end{thm}
\newpage

\begin{defn}Let $G$ be a group, and let $I$ and $S$ be non-empty sets, and for each $i \in I$ let $f_{i}:G \rightarrow S$ be a class function on $G$ (that is, if $g$ and $h$ are conjugate in $G$, then $f_{i}(g)=f_{i}(h)$, so that $f_{i}$ is constant on conjugacy classes of $G$). We say that the set $\{ f_{i} : i \in I \}$ of functions $ \textit{separates conjugacy classes}$ of $G$ if whenever $g$ and $h$ are in distinct conjugacy classes of $G$, then there exists an $i$ such that $f_{i}(g) \neq f_{i}(h)$.
\end{defn}

To prove Theorem 1.1, we will make use of the following theorem (see [3]):

\begin{thm} Let $G$ be a finite group, and let $S=\{g_{1}, g_{2}, g_{3}, ..., g_{k}\}$ be a complete set of representatives of the conjugacy classes of $G$. Then the following are equivalent: \newline
\newline
(1) $G$ is a rational group.\newline
\newline
(2) The permutation characters $\{1_{<g_{i}>}^{G}: i=1,2,...,k\}$ separate conjugacy classes \newline 
\indent $\hskip 3pt$ of $G$. \newline
\newline
(3) $G$ has a finite collection $\mathcal{H}$ of subgroups  such that the permutation characters \newline 
\indent $\hskip 3pt$ $\{1_{H}^{G}: H \in \mathcal{H}\}$ separate conjugacy classes of $G$.\newline
\newline
(4) $G$ has a permutation character which separates the conjugacy classes of $G$. \newline
\indent $\hskip 3pt$ (That is, G has a permutation representation such that any two elements of $G$ \newline
\indent $\hskip 3pt$  fixing the same number of letters are conjugate in $G$.)

\end{thm}

We can now give a proof of the main Theorem 1.1:

\begin{proof} (of Theorem 1.1) We will first show that (1) implies (2). Assume that $G$ is a rational group. Let  $\{g_{1}, g_{2}, g_{3}, ..., g_{k}\}$ be a complete set of representatives of the conjugacy classes of $G$, and let $I=\{1,2,3,...,k\}$. For each integer $i \in I$, let $X_{i}$ be the set of left cosets of $\langle g_{i} \rangle$  in $G$. Let $X = \bigcup_{i \in I} X_{i}$. (Without loss of generality, we will assume that the sets $X_{i}$ are pairwise disjoint, for otherwise we could replace $X_{i}$ with $X_{i} \times \{i\}$.) Define the type of each element of $X_{i}$ to be $i$. Define the incidence relation * on $X$ by $x*x'$ iff the cosets $x$ and $x'$ are equal or have different types and non-empty intersection.  This defines a finite geometry $\Gamma = (X, *, t)$. (In fact, $\Gamma$ is the coset geometry of the cyclic subgroups $\langle g_{i} \rangle$ for $i \in I$.) $G$ acts by left multiplication on the set $X$ of left cosets. So each $g \in G$ defines a function $\alpha_{g} : X \rightarrow X$ given by $\alpha_{g}(h\langle g_{i} \rangle) = gh\langle g_{i} \rangle$. This function $\alpha$ is a bijection with inverse $\alpha^{-1}=\alpha_{g^{-1}}$, and both $\alpha$ and $\alpha^{-1}$ preserve incidence and types. Therefore, $\alpha$ is an automorphism of $\Gamma$ and the mapping $g \mapsto \alpha_{g}$ defines a homomorphism $\alpha : G \rightarrow Aut(\Gamma)$; thus $G$ acts on $\Gamma$. Now assume that $g$ and $h$ are elements of $G$ fixing the same number of flags of type $J$ for all subsets $J$ of $I$ (and so in particular, for all $i \in I$,  the number of flags of type $i$ fixed by $g$  equals the number of flags of type $i$ fixed by $h$). This means that for each $i \in I$,  the number of cosets of $\langle g_{i} \rangle $ fixed by $g$ equals the number of cosets of $\langle g_{i} \rangle$ fixed by $h$. But $1_{\langle g_{i} \rangle }^{G}$ is the permutation character of $G$ acting on the cosets of $\langle g_{i} \rangle $, and so we have $1_{\langle g_{i} \rangle }^{G}(g) = 1_{\langle g_{i} \rangle }^{G}(h)$ for all $i \in I$. But by Theorem 1.3 above, since $G$ is a rational group, the permutation characters $\{1_{\langle g_{i} \rangle }^{G}: i \in I\}$ separate conjugacy classes of $G$. So $g$ must be conjugate to $h$ in $G$. So we have shown that (1) implies (2).\\

Next, we will show that (2) implies (1). Assume that there exists a finite geometry $\Gamma$ of type $I$ and action of $G$ on $\Gamma$ such that if $g$ and $h$ are elements of $G$ fixing the same number of flags of type $J$ for all subsets $J$ of $I$, then $g$ and $h$ are conjugate in $G$. We want to show that $G$ is a rational group. Let  $C =\{g_{1}, g_{2}, g_{3}, ..., g_{k}\}$ be a complete set of representatives of the conjugacy classes of $G$. Let $g_{i}$ and $g_{j}$ be distinct elements of $C$. Since $g_{i}$ and $g_{j}$ are not conjugate, there exists a subset $J$ of $I$ such that $g_{i}$ and $g_{j}$ fix a different number of flags of type $J$. It follows that there must exist a G-orbit $\mathcal{O}$ of flags of type $J$ on which $g_{i}$ and $g_{j}$ fix a different number of flags. Since $G$ is transitive on the G-orbit $\mathcal{O}$, the permutation character of $G$ on the G-orbit $\mathcal{O}$ is given by $1_{H}^{G}$, where $H$ is the stabilizer in $G$ of a flag of type $J$ in $\mathcal{O}$. So we have found a subgroup $H = H_{\{g_{i},g_{j}\}}$ (depending on the set  $\{g_{i}, g_{j}\}$) and a permutation character $1_{H}^{G}$ such that $1_{H}^{G}(g_{i}) \neq 1_{H}^{G}(g_{j})$. Repeating this process for all sets $\{g_{i}, g_{j}\}$ of distinct elements of $C$ gives us a finite collection $\mathcal{H}$ of subgroups such that the permutation characters $\{1_{H}^{G}: H \in \mathcal{H}\}$ separate conjugacy classes of $G$. By Theorem 1.3 above, it follows that $G$ is a rational group. Therefore (2) implies (1). This proves Theorem 1.1.
\end{proof}

Theorem 1.1 can be used to give a geometric proof of the well-known fact that the finite Coxeter groups of type $A_{n}$ (which are finite symmetric groups $Sym(\Omega)$) are rational groups. To do that, we will make use of the following lemma. Here, $Fix_{k}(g)$ represents the collection of subsets of $\Omega$ of cardinality $k$ fixed by the element $g$ of  $Sym(\Omega)$.

\begin{lem} 
Let $\Omega$ be a finite non-empty set, and let $g$ and $h$ be in $Sym(\Omega)$. If $|Fix_{k}(g)|= |Fix_{k}(h)|$ for every $k$, $(0 \leq k \leq |\Omega|)$, then $g$ and $h$ are conjugate in $Sym(\Omega)$.
\end{lem}
\noindent In other words, if two elements of $Sym(\Omega)$ fix the same number of subsets of $\Omega$ of cardinality $k$ for every $k$, $(0 \leq k \leq |\Omega|)$, then they are conjugate in $Sym(\Omega)$.

\begin{proof}
We will prove the equivalent contrapositive. Let $g$ and $h$ be in $Sym(\Omega)$ and assume that they are not conjugate in $Sym(\Omega)$. Then as permutations of $\Omega$, $g$ and $h$ do not have the same cycle structure. So there exists a least positive integer $m$ in $\{1,2,3,...,|\Omega|\}$ such that $g$ and $h$ have a different number of $m$-cycles. The number of  subsets of $\Omega$ of cardinality $m$ fixed by a permutation of $\Omega$ equals the number of ways of forming a set of cardinality $m$ from cycles of that permutation, which is different for  $g$ and $h$. Thus $|Fix_{m}(g)| \neq |Fix_{m}(h)|$. So it is not the case that $|Fix_{k}(g)|= |Fix_{k}(h)|$ for every $k$, $(0 \leq k \leq |\Omega|)$.
\end{proof}

\noindent \textit{Example:} Let $\Omega=\{1,2,3,4\}$, so $Sym(\Omega)$ acts on the 16 subsets of $\Omega$. The involution $(12)(34)$ fixes one subset of cardinality $0$ (the empty set), no subsets of cardinality $1$, two subsets of cardinality $2$ (namely $\{1,2\}$ and $\{3,4\}$), no subsets of cardinality $3$, and one subset of cardinality $4$ (namely $\Omega$ itself). Also, the involution $(13)(24)$ fixes one subset of cardinality $0$, no subsets of cardinality $1$, two subsets of cardinality $2$, no subsets of cardinality $3$, and one subset of cardinality $4$. So the lemma tells us that $(12)(34)$ and $(13)(24)$ are conjugate in $Sym(\Omega)$. On the other hand, the permutations $(12)(34)$ and $(1234)$ are not conjugate in $Sym(\Omega)$, and so by the lemma, they must fix different numbers of sets of cardinality $k$, for some $k$. Checking, we see that they both fix the same number of sets of cardinalities $0, 1, 3,$ and $4$, but while $(12)(34)$ fixes two subsets of cardinality 2, the permutation $(1234)$ fixes no subset of cardinality $2$. Note also that the elements $(12)(34)$ and $(123)$ each fix four subsets of $\Omega$, but are not conjugate; so it is not possible to weaken the hypotheses of Lemma 1.4 to the condition that $g$ and $h$ only fix the same total number of subsets of $\Omega$.\\

To give a geometric proof of the well-known fact that the finite Coxeter groups of type $A_{n}$ (which are finite symmetric groups $Sym(\Omega)$) are rational groups, we will construct a finite geometry $\Gamma$ on which $Sym(\Omega)$ acts, and then apply Theorem 1.1. Let $\Omega = \{1,2,...,n \}$ and let $X = \mathcal{P}(\Omega)$, the power set of $\Omega$. For each element $x$ of $X$ (i.e. for each subset $x$ of $\Omega$), define the type $t(x)$ of $x$ to be the cardinality of $x$. For elements $x$ and $x'$ in $X$, define $x*x'$ if and only if either $x \subseteq x'$ or $x' \subseteq x$. Then $\Gamma =(X, *, t)$ is a finite geometry of type $I$, where $I=\{0,1,2,...,n \}$. The usual action of $Sym(\Omega)$ on $\Omega$ induces a natural action of $Sym(\Omega)$ on the geometry $\Gamma$.  We claim that this action satisfies condition (2) of the main Theorem 1.1: let $g$ and $h$ be arbitrary elements of $Sym(\Omega)$ which fix the same number of flags of type $J$ for all subsets $J$ of $I$; we want to show that $g$ and $h$ are conjugate. For each  integer $i \in  I=\{0,1,2,...,n \}$,  the elements $g$ and $h$ fix the same number of flags of $\Gamma$ of type $i$, and therefore fix the same number of subsets of $\Omega$ of cardinality $i$. Then by  Lemma 1.4, $g$ and $h$ are conjugate in $Sym(\Omega)$. So this action satisfies condition (2) of the main Theorem 1.1, and so $Sym(\Omega)$ is a rational group.\\

\noindent \textit{Concluding question}: For each of the other previously mentioned finite Coxeter groups, $B_{n}, D_{n}, E_{6}, E_{7}, E_{8}, F_{4},$ and $G_{2}$, can one - without assuming rationality - find a natural geometry (perhaps related to its geometry of maximal parabolic subgroups) which satisfies the geometric condition of Theorem 1.1? If so, this would provide a natural geometric proof of the rationality of these groups.

\bibliographystyle{amsplain}

\end{document}